\documentclass[onefignum,onetabnum]{siamart171218}

\usepackage{lipsum}
\usepackage{amsfonts}
\usepackage{graphicx}
\usepackage{epstopdf}
\usepackage{algorithm} 
\usepackage[noend]{algpseudocode} 
\usepackage{psfrag}
\usepackage{pstool}

\usepackage{float} 
\usepackage{placeins} 
\usepackage{booktabs} 
\usepackage{subfig} 

\ifpdf
  \DeclareGraphicsExtensions{.eps,.pdf,.png,.jpg}
\else
  \DeclareGraphicsExtensions{.eps}
\fi

\newsiamremark{remark}{Remark}
\newsiamremark{hypothesis}{Hypothesis}
\crefname{hypothesis}{Hypothesis}{Hypotheses}
\newsiamthm{claim}{Claim}

\headers{Efficient p-multigrid spectral element model }{A. P. Engsig-Karup and W. L. Laskowski}

\title{Efficient p-multigrid spectral element model for water waves and marine offshore structures\thanks{Submitted to the editors August 24th, 2020.
}}

\author{Allan P. Engsig-Karup\thanks{Technical University of Denmark, 2800 Kgs. Lyngby, Denmark
  (\email{apek@dtu.dk}, \url{http://www.compute.dtu.dk/\string~apek/}).}
\and Wojciech Laskowski\thanks{Department of Applied Mathematics and Computer Science, Technical University of Denmark, (\email{apek@dtu.dk}).}
}

\usepackage{amsopn}

\ifpdf
\hypersetup{
  pdftitle={An efficient p-multigrid spectral element model for fully nonlinear water waves and fixed bodies},
  pdfauthor={A. P. Engsig-Karup, and W. L. Laskowski}
}
\fi

\begin{document}

\maketitle

\begin{abstract}
In marine offshore engineering, cost-efficient simulation of unsteady water waves and their nonlinear interaction with bodies are important to address a broad range of engineering applications at increasing fidelity and scale. We consider a fully nonlinear potential flow (FNPF) model discretized using a Galerkin spectral element method to serve as a basis for handling both wave propagation and wave-body interaction with high computational efficiency within a single modelling approach. We design and propose an efficient $\mathcal{O}(n)$-scalable computational procedure based on geometric $p$-multigrid for solving the Laplace problem in the numerical scheme. The fluid volume and the geometric features of complex bodies is represented accurately using high-order polynomial basis functions and unstructured meshes with curvilinear prism elements. The new $p$-multigrid spectral element model can take advantage of the high-order polynomial basis and thereby avoid generating a hierarchy of geometric meshes with changing number of elements as required in geometric $h$-multigrid approaches. 
 We provide numerical benchmarks for the algorithmic and numerical efficiency of the iterative geometric $p$-multigrid solver. Results of numerical experiments are presented for wave propagation and for wave-body interaction in an advanced case for focusing design waves interacting with a FPSO. Our study shows, that the use of iterative geometric $p$-multigrid methods for the Laplace problem can significantly improve run-time efficiency of FNPF simulators.
\end{abstract}

\begin{keywords}
High-order numerical method, spectral element method, fully nonlinear potential flow, Laplace problem, geometric $p$-multigrid, marine offshore hydrodynamics.
\end{keywords}

\begin{AMS}
  65B99, 65M55, 76B07 
\end{AMS}

\section{Introduction}

The description of water waves and their nonlinear interaction with bodies are of practical importance in marine offshore hydrodynamics. Few existing  phase-resolving simulators are capable of handling both the wave propagation and wave-body interaction problem in an computationally efficient and scalable way, and within the same  simulator. Adding such capability is, e.g., of relevance in many applications within offshore engineering and increasingly in the renewable sectors. For example, for designing offshore wind turbine foundations \cite{BREDMOSE2016379}, and improving wave energy device simulations \cite{DavidsonCostello2020}. Research in theoretical foundations of wave-body modelling \cite{John1949,lannes2016dynamics}, improved uncertainty  analysis approaches \cite{BEE2016,ehi2019activesubspace} and practical applications \cite{BOSI2019222} illuminates the  opportunities and demonstrate significant potential for improving the fidelity of conventional practical engineering tools for engineering analysis using nonlinear wave models. In water wave simulators, describing the fluid volume and taking into account bodies or structures in the representation can be handled using unstructured meshes and use high-order methods to improve cost-efficiency of simulations and reach larger scales to better estimate the influence of sea states in coastal areas. For water wave simulations, there are still progress to be made in improving dispersive and nonlinear wave propagation models, e.g to use unstructured numerical methods, e.g. see review \cite{Brocchini2013}. Also, recent research shows progress in the development of high-order accurate methods for nonlinear and dispersive large-scale ocean simulations  \cite{engsigkarupetal2012,engsigkarup2017stabilised,BrusEtAl2019,GlimbergEtAl2019}. 

For unsteady water wave modelling, the fully nonlinear potential flow (FNPF) model represents one of the two most widely used classes of nonlinear models for water waves and is particularly important due to the ability to solve the water wave propagation problem while accounting for dispersion and nonlinear transformations of waves. This can be done with much higher efficiency than the higher fidelity CFD solvers that is based on Navier-Stokes type models when the rotational and viscus effects of the fluid motion are negligible \cite{RansleyEtAl2018}. The most cost-efficient numerical schemes for time-dependent problems with long integration times are those that are of high-order accuracy \cite{KO72}. Thus, numerical schemes that have fast  convergence rates and scale linearly with problem size (measured in terms of degrees of freedom in the discretization) has significant potential to reduce the cost of CFD and enable simulations of increasing fidelity in practical times  \cite{DavidsonCostello2020}. A proper designed multigrid scheme for the linear solver in CFD tools is a key component for scalability and high efficiency \cite{HuismanStillerJochen2018}. Furthermore, the majority of practical simulation scenarios requires adaptive unstructured grids. Other high-order numerical methods that addresses the wave propagation problem and the wave-body problem in a single solver strategy is the high-order boundary element method (HOBEM) \cite{Ferrant1996,HagueSwan2009,HarrisEtAl2014}
that is particularly strong in handling the geometry using unstructured grids, but is limited in terms of numerical efficiency due to inefficient asymptotic scaling of work effort as a result of high computational complexity in the discrete solution of the resulting system of equations in the solver. The SWENSE technique alleviates the efficiency problem of CFD tools through a decomposition approach  \cite{DucrozetEtAl2014,ZhaobinEtAl2018}. However, it is less suitable for wave propagation covering long duration of times or large domains due to lack of incident wave solutions for the wave propagation problem in varying bathymetry and that numerical dispersion errors grows over time. Furthermore, the use of high-order numerical methods for dispersive and nonlinear free surface hydrodynamics are still surprisingly scarce, especially when structural bodies are taken into account. Two widely used FNPF solvers for nonlinear wave propagation are the high-order finite differences method (FDM) that is both efficient and scalable   \cite{ENGSIGKARUP20092100,EngsigKarup2014,GlimbergEtAl2019} and the high-order spectral method (HOS) \cite{WestEtAl1987,DOMYOU1987}. The FDM discretizations can be tailored to use adaptive and body conforming meshes, e.g. via multiblock discretization and curvilinear grids \cite{GlimbergEtAl2019}, methods such as HOS based on approximations via global basis functions do not have this flexibility.

In this work, we seek to address the computational bottleneck problem of FNPF solvers. We propose a geometric $p$-multigrid ($p$-GMG) accelerated iterative Laplace equation solver that is a fundamentally important extension of a recently developed spectral element solver that have been successfully validated against experimental measurements \cite{EngsigKarupEskilsson2018,RansleyEtAl2019,EngsigKarupEskilsson2019}. The solver strategy is based on the spectral element method attributed to Patera \cite{PAT84} suitable for efficient and scalable simulations for elliptic problems  \cite{PavarinoWarburton2000}. A recent review of state-of-the-art spectral element methods is given in \cite{XuiEtAl2018}. The single-phase FNPF model is based on explicit tracking of the water surface, and avoids the difficulties associated with interface tracking solvers such as volume of fluid (VOF) \cite{Ubbink1997,Jasak07openfoam:a,Refresco2009} and level set (LS) method \cite{GroossHesthaven2005,BIHS2016191,KarakusEtAl2018}. VOF methods are among the most widely used general-purpose solvers for engineering applications, while LS methods \cite{OSHER198812,BIHS2016191} sees increasing maturity for free surface flows. These types of solvers comes with the limitation that it is difficult to device high-order accurate numerical schemes that captures the complex interfaces that is needed to make them accurate and cost-efficient. This makes them prone to errors in conserved properties such as mass and energy,  significant numerical damping in wave-propagation problems, which historically have limited these solvers to applications in relatively small domains and relatively short integration times for accurate prediction of wave propagation and kinematics. Also, due to the low-order accuracy of the methods, significant number of cells is needed per wave length of the highest harmonics that needs to be resolved to deliver converged results for engineering purposes. Still, these CFD solvers are widely accepted to be favourable for the wave-body problem due to the possibility for describing bodies of arbitrary geometric  complexity and accounting for violent effects such as breaking waves, slamming, aeration effects (multi-phase), green water on ship decks, etc. that is needed in typical industrial applications for estimation of wave-induced loads. 

The spectral element methods (e.g. see \cite{PAT84,KS99,deville_fischer_mund_2002,XuiEtAl2018}) is a class of finite element methods with support for adaptive meshes and that is based on piece-wise continuous arbitrary order polynomial approximations for spatial discretization. These methods are attractive for a broad range of applications due to their exponential convergence rates and geometric flexibility. In the lat decade, the spectral element methods have seen increasing maturity and adoption in terms of real-world applications. However, the main limitations are still in the mesh generation \cite{TURNER2016340}, and the work effort associated with the sparse operators with dense blocks for each element that calls for use of matrix-free implementations and high-performance computing \cite{MarkallEtAl2013}. Remark, that multigrid methods can be accelerated via mixed-precision strategies on modern many-core architectures \cite{GlimbergEtAl2011,engsigkarupetal2012,abdelfattah2020survey} to gain additional performance in parallel implementations.

In the context of FNPF solvers, recent work have demonstrated for practical applications that the methodology can be stabilised and is accurate \cite{EEB2016,EngsigKarupEskilsson2018}. A remaining scientific challenge for enabling large-scale 3D simulations for FNPF solvers is the development of an efficient iterative solver strategy that is efficient for wave-body problems. A proper designed multigrid solver strategy can achieve the ideal $\mathcal{O}(n)$ scaling of work effort, cf. \cite{BRANDT77,Brandt1998,Trottenberg01}, and the most efficient multigrid solvers achieve "textbook multigrid efficiency" (TME), which is defined by Brandt \cite{Brandt1998} as \textit{solving a discrete PDE problem in a computational work which is only a small (less than 10) multiple of the operation count in the discretized system of equations itself}. This is usually measured in terms of residual calculations for the discretized system. The work effort to achieve satisfactory approximate numerical solution depends on both the numerical efficiency of the multigrid schemes and the tolerance chosen. In solving linear systems of equations using multigrid methods, to maximize attainable accuracy (cf. \cite{EngsigKarup2014}) one need to reduce the algebraic errors below truncation errors, since from the triangle inequality for the error between the true solution ($\phi$) and the approximate solution ($\hat{\phi}_h$) to a system of equation is given as
\begin{align}
    \underbrace{||\phi-\phi_h||}_\text{Errors}\leq \underbrace{|| \phi-\hat{\phi}_h ||}_\text{Truncation errors} + \underbrace{|| \hat{\phi}_h - \phi_h ||}_\text{Algebraic errors},
\end{align}
and hence through the use of iterative solvers the attainable accuracy is limited by the level of truncation errors with the level of error influenced by the algebraic errors ($\hat{\phi}$). Thus, the scaling properties and TME set a clear optimality target to define to what extend a given solver strategy is efficient. The main objective of this work has been to design a geometric $p$-multigrid accelerated FNPF solver strategy and clarify to what extend TME can be met. 

There are two classes of multi-level methods that form the basis for multigrid schemes, namely, the algebraic multigrid (AMG) and geometric multigrid (GMG) methods. The AMG methods are among the most widely used in industrial applications as it handle the solution of system of equations as a black box iterative solver strategy \cite{BrandtEtAl1984,DENDY198657,RugeStuben1987,OlSc2018} with the fine grid discretization of the model problem as input. AMG then proceeds by constructing coarser representation algebraically according to certain rules. However, when dealing with unsteady numerical solvers where grids may change and as a result the system of equations or its coefficients will be subject to changes, then AMG falls short as an efficient time-domain solver strategy due to the need for updates of the operators which happens in a post-processing step before the solve procedure. In contrast, the GMG method starts from a coarse grid problem and then refine the mesh resolution through a mesh hierarchy according to certain rules. GMG can be designed to accommodate updates of local element matrices in a more straightforward way, utilize matrix-free operator implementations, and as a result achieve a low computational cost for unsteady problems that requires frequent operator updates.
 


Multigrid methods can be used to accelerate the convergence rate of iterative schemes by exploiting the {\em smoothing property} where oscillatory modes of the error are damped fast and smooth modes are damped slowly. This property can be taken advantage of in iterative schemes by using a nested grid consisting of fine and coarser grids. An iterative scheme having the smoothing property is called a {\em smoother}. Any multigrid algorithm involves components for smoothing, restriction, interpolation, solution on different grids and choice of cycle type. 
The framework for multigrid methods is general, however, usually the components of multigrid methods have to be tailored to specific problems \cite{Brandt1984}. Multigrid algorithms are attractive as they provide the basis for solving large system of equations with complexity given as $\mathcal{O}(n)$ and make it possible to achieve scalable solution effort, i.e. effort which only increases linearly with the degrees of freedom $n$ in the discretization. 







%


The paper is organised as follows. The governing equations for fully nonlinear potential flow are described in \cref{sec:gov}. \Cref{sec:num} describes the key ideas behind a geometric $p$-multigrid model that is proposed. In \cref{sec:itsol} numerical results from benchmarks of the new  iterative $p$-GMG solver is provided and \cref{sec:alg}  contains results from the application of the geometric $p$-multigrid solver in benchmarks in 2D and 3D using two different FNPF model formulations. Results are given for the assessment of the performance in terms of efficiency and scalability achieved for different $p$-GMG solver strategies for the practical benchmarks. The conclusions, follow in \cref{sec:conclusions}.


\section{Governing Equations}
\label{sec:gov}

Irrotational, inviscid fluid flow and non-breaking water waves are governed by the Fully Nonlinear Potential Flow (FNPF) equations defined on a fluid domain $\Omega$. Domain boundaries may be constituted of a varying bathymetry ($\Gamma^b$), the free surface ($\Gamma^{FS}$), horizontal boundaries, e.g. in the form of walls of a numerical wave tank, and internal boundaries in the form of the faces of structural bodies or fixed structures. The non-dimensional FNPF model equations subject to kinematic and dynamic free surface conditions and the continuity equation on $\Omega$ can be written in an Eulerian formulation \cite{EGNL13} as 
\begin{subequations}
\begin{align}
\frac{\partial \eta}{\partial t}  &= -\boldsymbol{\nabla}\eta\cdot\boldsymbol{\nabla}\tilde{\phi}+\tilde{w}(1+\boldsymbol{\nabla}\eta\cdot\boldsymbol{\nabla}\eta) \quad \quad \quad \quad \quad \quad \textrm{in} \quad \Gamma^{\textrm{FS}}\times T,  \label{FSeta} \\
\frac{\partial \tilde{\phi}}{\partial t} &= -g\eta - \frac{1}{2}\left(\boldsymbol{\nabla}\tilde{\phi}\cdot\boldsymbol{\nabla}\tilde{\phi}-\tilde{w}^2(1+\boldsymbol{\nabla}\eta\cdot\boldsymbol{\nabla}\eta)\right) \quad \textrm{in} \quad \Gamma^{\textrm{FS}} \times T,
\label{FSphi}
\end{align}
\end{subequations}
where the '$\sim$' symbol is used for free surface variables. Closure for the free surface equations are obtained by solving the Laplace problem
\begin{subequations}
\begin{align}
\phi &=\tilde{\phi}, \quad z=\eta \quad \textrm{on}\quad \Gamma^{FS}, \\
\nabla^2 \phi &= 0, \quad -h(x,y)<z<\eta \quad \textrm{in} \quad \Omega, \\
{\bf n}\cdot \nabla\phi &= 0, \quad z=-h({x,y})\quad \textrm{on}\quad \Gamma^b.
\end{align}
\label{eq:laplaceproblem}
\end{subequations}
In our numerical benchmarks we consider both the fully nonlinear free surface formulation \eqref{FSphi} and the small-amplitude version ($H/L\approx0$) of this free surface problem. Here $L$ is the wave length and $H$ is the wave amplitude. The gravitational acceleration $(g)$ is assumed to be 9.81 $m/s^2$.



\section{Numerical method}
\label{sec:num}

\subsection{Spectral Element Method}

For the space discretization, we use the implementation of a spectral element method detailed in 2D for a $\sigma$-transformed version \cite{EEB2016} and in 3D for a non $\sigma$-transformed version \cite{EngsigKarupEskilsson2018}  for the governing equations (see Section \ref{sec:gov}). These already validated solvers serve the purpose of benchmarking the geometric $p$-multigrid solver in this study for solving the Laplace problem. We revisit the discretization details in the following.

The Laplace problem can be formulated in the general non $\sigma$-transformed form most suitable for wave-body applications via directly discretizing \eqref{eq:laplaceproblem}
or alternatively via the $\sigma$-transformed form that is most suitable for wave propagation applications 
\begin{subequations}
\begin{align}
\phi &=\tilde{\phi}, \quad z=\eta \quad \textrm{on}\quad \Gamma^{FS}, \\
\nabla^\sigma\cdot(\mathcal{K}({\bf x;t})\nabla^\sigma\phi) &= 0, \quad {\bf x}\in\Omega^\sigma, \\
{\bf n}\cdot \nabla\phi &= 0, \quad z=-h({x,y})\quad \textrm{on}\quad \Gamma^b,
\end{align}
\end{subequations}
where $\nabla^\sigma=(\partial_x,\partial_y,\partial_\sigma)$ is the differential operator that act in the reference domain with a $\sigma$-coordinate.

The weak formulation of \eqref{eq:laplaceproblem} is stated in the form
\begin{align}
\int_\Omega \nabla^2 \phi d{\bf x} = \oint_{\Gamma} v {\bf n} \cdot \nabla \phi d{\bf x} - \int_\Omega \nabla \phi \cdot \nabla v d{\bf x}.
\label{eq:laplaceeq}
\end{align}
Let $\Omega^c\subset \mathbb{R}^3$ be the time-independent computational domain $\Omega^c=\{ (\boldsymbol{x},\sigma) | \boldsymbol{x}\in\Gamma^{\textrm{FS}}, 0\leq \sigma \leq 1 \}$. If we introduce the Jacobian of the map $\Psi:\Omega\to\Omega^{c}$ 
\begin{align}
\label{eq:jacobian}
\mathcal{J}(\boldsymbol{x},z,t) =  \left[  
\begin{array}{ccc} 
\partial_x x & \partial_y x  & \partial_z x \\
\partial_x y & \partial_y y  & \partial_z y \\
\partial_x \sigma & \partial_y \sigma  & \partial_z \sigma
\end{array} 
\right]
=
 \left[  
\begin{array}{ccc} 
1 & 0 & 0 \\ 
0 & 1 &  0 \\ 
\frac{\partial_xh}{d}-\frac{\sigma \partial_xd}{d}  & \frac{\partial_yh}{d}- \frac{\sigma \partial_yd}{d} & \frac{1}{d} 
\end{array} 
\right],
\end{align}
then the $\sigma$-transformed formulation \cite{EEB2016} that is appropriate for solving efficiently the wave propagation problem in a FNPF model can be stated in the form
\begin{align}
\mathcal{K}(\boldsymbol{x};t) = \frac{1}{\textrm{det} \mathcal{J}}\mathcal{J J}^T =  \left[  
\begin{array}{ccc}
d & 0 & -\sigma \frac{\partial\eta}{\partial x}\\ 
0 & d & -\sigma \frac{\partial\eta}{\partial y}\\ 
-\sigma \frac{\partial\eta}{\partial x}  & -\sigma \frac{\partial\eta}{\partial y} & \frac{1+(\sigma \frac{\partial\eta}{\partial x})^2+(\sigma \frac{\partial\eta}{\partial y})^2}{d}
\end{array} \right],
\end{align}
or in the non $\sigma$-transformed version \eqref{eq:laplaceproblem} that is suited for wave-body applications. In both cases, the resulting discretised systems are derived from a weak formulation of the form
\begin{align}
\int_{\Omega^\sigma} \nabla^\sigma \cdot (\mathcal{K} \nabla^\sigma)\phi v d{\bf x} = \oint_{\Gamma} v {\bf n}\cdot(\mathcal{K} \nabla^\sigma) \phi d{\bf x}- \int_{\Omega^\sigma} \mathcal{K} \nabla^\sigma \phi v d{\bf x} = 0.
\end{align}

The details of the spatial Galerkin discretization leading to a Spectral Element Method is briefly described next. We start by forming a partition of the spatial domain $\Omega$ by a tessellation $\mathcal{T}_h$ of the $xy$-plane into $N^k$ non-overlapping shape-regular elements $\mathcal{T}_k$ such that $\cup_{k=1}^{N^{k}}\mathcal{T}_k=\mathcal{T}_h$ with $k$ denoting the $k$'th element. We then introduce the space of continuous, piece-wise polynomial functions 
\begin{align}
V^p=\{ v_h\in C^0(\Omega); \forall k \in \{ 1,...,N^k \}, v_{h|\mathcal{T}_k}\in\mathbb{P}^p \},
\end{align}
where $\mathbb{P}^p$ is the space of polynomials of degree at most $p$. In the Galerkin scheme, we choose test functions $v\in V^p$.

We introduce the finite-dimensional approximations
\begin{align}
\phi_h^p(\boldsymbol{x},t) = \sum_{i=1}^{N^{p}} \phi_i^p(t) N_i^p(\boldsymbol{x}),
\end{align}
where $\{N_i^p(\boldsymbol{x})\}_{i=1}^{N^{p}}\in V^p$ is the finite set of global finite element basis functions with cardinal property $N_i^p({\boldsymbol{x}_j^p})=\delta_{ij}$ at mesh nodes $\boldsymbol{x}_j^p\in\Omega^p$ with $\delta_{ij}$ the Kronecker Symbol. If we substitute these expressions into the weak formulation and choose $v^p(\boldsymbol{x})\in\{ N_i^p(\boldsymbol{x}) \}_{i=1}^{N^{p}}$ the discretization can be made symmetric and results in a sparse linear system of the form
\begin{subequations}
\begin{align}
\mathcal{L}^p \phi_h^p = {\bf b}^p, \quad \mathcal{L}^p\in\mathbb{R}^{N^p\times N^p}, \quad \phi_h^p, {\bf b}^p\in\mathbb{R}^{N^p},
\label{eq:laplaceproblemdiscrete}
\end{align}
where $N^p$ is the total degrees of freedom in the discretisation in the $p$'th mesh and
%
\begin{align}
\mathcal{L}_{ij}^p &= -\int_{\Omega^c} (\mathcal{K} \nabla^c N_j^p)  \cdot \nabla^c N_i^p d\boldsymbol{x} = -\sum_{k=1}^{N^k}   \int_{{\Omega^{c,k}}} (\mathcal{K} \nabla^c N_j^p)  \cdot \nabla^c N_i^p d\boldsymbol{x}, \\
b_i^p &= \oint_{\Gamma^c} N_i^p {\bf n}\cdot \left(\mathcal{K}\nabla^c\left( \sum_{j=1}^{N^p} \phi_j^p N_j^p \right)\right) d\boldsymbol{x} \nonumber \\
&= \sum_{s=1}^{N^s} \oint_{\Gamma^{c,s}} N_i^p {\bf n}\cdot \left(\mathcal{K}\nabla^c\left( \sum_{j=1}^{N^p} \phi_j^p N_j^p \right)\right) d\boldsymbol{x},
\end{align}
\end{subequations}
and where $N^s$ is total number of segments of the domain boundary of the computational domain constituted from the surfaces of the elements used to tesselate the domain.

If the preconditioned conjugate gradient (PCG) is used for the solution, convergence is guaranteed \cite{CaiEtAl1998}. Furthermore, the convergence properties are determined by the eigenspectrum $\rho(\mathcal{L})$ which for symmetric positive definite (SPD) systems are characterized by eigenvalues $\lambda\in[\lambda_{min};\lambda_{max}]$. Focus is given to the design of an efficient iterative $p$-multigrid spectral element solver for the Laplace problem  \eqref{eq:laplaceproblemdiscrete}.

\section{Iterative solver}
\label{sec:itsol}

The goal is to solve \eqref{eq:laplaceproblemdiscrete} using a $\mathcal{O}(n)$-scalable iterative $p$-multigrid iterative scheme tailored for high-order accurate approximations on fixed meshes where  $n$ is the total degrees on the finest mesh. An essential component for those methods is to have a coarse grid acceleration, which results in $p$-multigrid strategies \cite{FL05}. We solve the linear system of equations that arises from the SEM discretization, cf. \eqref{eq:laplaceproblemdiscrete}. 

\subsection{Preconditioning}

Preconditioning implies a transformation of the original problem \eqref{eq:laplaceproblemdiscrete}
%
%
to a new problem with same solution, e.g. of the left-preconditioned form
\begin{align}
\mathcal{M}^{-1}\mathcal{L}\phi_h=\mathcal{M}^{-1}{\bf b}, \quad \mathcal{M}\in{\bf R}^{N^p\times N^p},
\label{eq:precondsys}
\end{align}
where $\mathcal{M}$ is the preconditioner. Identifying a proper preconditioning strategy is the key to efficient solution of sparse linear systems such as those that arise from numerical disceretization using the spectral element methods, e.g. discretization of \eqref{eq:laplaceproblemdiscrete}. 
The objective of preconditioning \cite{vanderVorst:2003:IKM} is to
is to accelerate the convergence of the iterative solver to make it less expensive to solve the system of equations.
In this work, the design of the preconditioner comes with the  requirement, that it should work efficiently in the setting of unstructured meshes and arbitrarily shaped bodies. For classical finite element methods (FEM), the main path to convergence is $h$-type convergence where the mesh is refined to reduce the spatial errors. In the spectral element method, an additional option is available, namely, $p$-type convergence, where the order of the local polynomial expansions to each element in the mesh is adaptive. The $p$-type convergence comes with the additional property that the spatial resolution can be improved by modification of the expansion basis and therefore do not need to change the initial mesh as in the FEM to improve on the errors. The objective of this study is to take advantage of this and then consider the design of a geometric $p$-multigrid scalable iterative solver strategy that works without $h$-refinement to avoid dealing with the details of the features of complex geometry of bodies represented via an adapted mesh. 

%

Previous benchmarks of near-optimal multigrid efficiency for a multigrid preconditioned defect correction method (PDC) obtained using a flexible-order finite difference FNPF simulator \cite{EngsigKarup2014} can be used as a reference for the multigrid methods in this work. In the following, we provide details on the design of an iterative solver strategy for the intended applications and demonstrate through numerical examples the efficiency and scalability of geometric $p$-multigrid solvers tailored to solving the FNPF equations in \cref{sec:gov}.
 
\section{Multigrid algorithm}
\label{sec:alg}

\begin{figure}[H]
\centering
\includegraphics[width=0.2\textwidth]{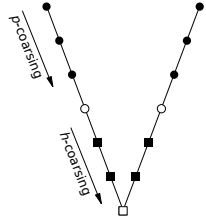}
\caption{Illustration of a multigrid $V$-cycle combining a $p$- and $h$-multigrid strategy.}
\label{fig:Vcycle}
\end{figure}

In this work we employ the standard version of the geometric  multigrid algorithm, namely the \textit{V-cycle} version  \cite{Trottenberg01}, i.e. starting from the finest level, coarsening until the coarsest grid is reached and going back to the finest level. This is depicted conceptually in Figure \ref{fig:Vcycle}, where a combination of $p$- and $h$-multigrid refinements is presented where each of the grids are visited using the recursive multigrid algorithm. Geometric multigrid methods relies on several grids ranging from the coarsest to the fine grids. The changes in resolution between grids are referred to coarsening. In standard coarsening, the coarser grids are simply half the size in each spatial direction of the finer grid. However, in some cases it can be advantageous to employ semi-coarsening where the grid is only coarsened in one spatial direction or combinations hereof.

\begin{figure}[H]
\centering
\includegraphics[width=0.5\textwidth]{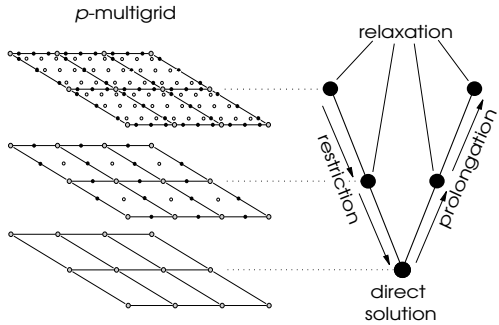}
\caption{Illustration of a $V$-cycle of geometric $p$-multigrid in 2D. Big, gray dots represent element's vertices (the $h$-mesh), black and white nodes represent Legendre Gauss Lobatto (LGL) nodes respectively on the edges and in the interior part of the element.}
\label{fig:pcoarse}
\end{figure}

We restrict our focus to the geometric $p$-multigrid version in the following. Consider the span of $P$ grids, discretized on the same domain $\Omega$ consisting of the same $N^k$ elements, however, utilising different polynomial order of the element expansions $N^p$, $p \in \{1,...,P\}$ which results in discretization pairs of the form $(N^p,N^k)$.  Typically, we strive to have coarsening between grids fulfilling $N^p \approx 2N^{p-1}$. 
This range of $p$-grids form the basis for a nested grid approach. For example, a $V$-cycle of a $p$-multigrid method on several grids $\Omega_p$ is depicted in figure \ref{fig:pcoarse}.

A pseudo-code for the recursive multigrid $V$-cycle for solving a linear system $\mathcal{A} {\bf x}={\bf f}$ is given as algorithm \ref{alg:mg}, employed with the following input: current grid level $p$, matrix $\mathcal{A}$ and right hand side ${\bf f}$ on a $p$ grid, initial approximation ${\bf u}_{h,p}$ and number of smoothing operations $\nu_1$ (pre-smoothing), $\nu_2$ (post-smoothing).

\begin{algorithm}
\caption{Geometric $p$-multigrid $V$-cycle}\label{alg:mg}
\begin{algorithmic}[1]
\Function{MultigridCycle}{$p,{\bf u}_{h,p},{\mathcal{A}}_p,{\bf f}_{h,p},{\nu}_1,{\nu}_2$}
\State ${\bf u}_{h,p} \gets {\bf u}_{h,p} + {\mathcal{S}}^{-1}({\bf f}_{h,p} - {\mathcal{A}_p} {\bf u}_{h,p}), \quad m=0,...,\nu_1-1$ 
\State ${\bf r}_{h,p} \gets {\bf f}_{h,p} - {\mathcal{A}}_p {u}_p$ 
\State ${\bf r}_{h,p-1} \gets {\mathcal{R}}_p^{p-1} {\bf r}_{h,p}$ 
\If{$p=1$}
\State $ \tilde{\bf u}_{h,p-1} \gets {({\mathcal{A}}_{p-1})}^{-1} {\bf r}_{h,p-1}$ 
\ElsIf{$p>1$}
\State $\tilde{\bf{u}}_{h,p-1} \gets \Call{MultigridCycle}{p-1,{0},{\mathcal{A}}_{p-1},{\bf r}_{h,p-1},{\nu}_1,{\nu}_2}$
\EndIf
\State $\tilde{{\bf u}}_{h,p} \gets {\mathcal{P}}_{p-1}^{p} \tilde{{\bf u}}_{h,p-1}$ 
\State ${\bf u}_{h,p}\gets {\bf u}_{h,p} + \tilde{{\bf u}}_{h,p}$ 
\State ${\bf u}_{h,p} \gets {\bf u}_{h,p} + {\mathcal{S}}^{-1}({\bf f}_{h,p} - {\mathcal{A}}_p {\bf u}_{h,p}), \quad m=0,...,\nu_2-1$ \\
\Return ${\bf u}_{h,p} \gets {\bf u}_{h,p}$
\EndFunction
\end{algorithmic}
\end{algorithm}

Multigrid iteration can be then combined with other iterative solvers in the form of preconditioned defect correction (PDC) \eqref{alg:mgdc} and preconditioned conjugate gradient (PCG) \eqref{alg:mgcg}. We use the stop criterion for the iterative solver in the form
%
$\|{\bf r}^{[m]}\| \leq  \texttt{rtol} \cdot {\|{\bf f}\|} + \texttt{atol}$, 
%
giving the user control to decide the acceptable level for convergence in terms of combined absolute (\texttt{atol}) and relative (\texttt{rtol}) error tolerances. 

\FloatBarrier
\begin{algorithm}
\caption{Multigrid Preconditioned Defect Correction}\label{alg:mgdc}
\begin{algorithmic}[1]
\Function{PDCpMG}{${\bf u}_h,\mathcal{A},{\bf f}_h,\texttt{atol},\texttt{rtol},i_{max},p_{grid},\nu_1,\nu_2$}
\While{$i < i_{max} \; \& \; {\bf r}_h > \texttt{rtol} \cdot {\|{\bf f}_h\|} + \texttt{atol} $}
\State ${\bf r}_h \gets {\bf f}_h - \mathcal{A} {\bf u}_h$
\State $\delta \gets \Call{MultigridCycle}{p_{grid},0,\mathcal{A},{\bf -r}_h,\nu_{1},\nu_{2}}$
\State ${\bf u}_h \gets  {\bf u}_h - \delta$ 
\State $i \gets i + 1$
\EndWhile \\
\Return $u$
\EndFunction
\end{algorithmic}
\end{algorithm}
\FloatBarrier

\FloatBarrier
\begin{algorithm}
\caption{$p$-Multigrid Preconditioned Conjugate Gradient}\label{alg:mgcg}
\begin{algorithmic}[1]
\Function{PCGpMG}{${\bf u}_h,\mathcal{A},{\bf f}_h,\texttt{atol},\texttt{rtol},i_{max},p_{grid},\nu_1,\nu_2$}
\State ${\bf r}_h \gets {\bf f}_h - \mathcal{A}{\bf u}_h$
\State $p \gets \Call{MultigridCycle}{p_{grid},0,\mathcal{A},{\bf r}_h,\nu_{1},\nu_{2}}$
\State $\delta \gets p^T {\bf r}_h$
\While{$i < i_{max} \; \& \; {\bf r}_h > \texttt{rtol} \cdot {\|{\bf f}_h\|} + \texttt{atol} $}
\State $q \gets \mathcal{A}{\bf \delta}$
\State $\alpha \gets \delta / (p^T q)$
\State ${\bf u}_h \gets {\bf u}_h + \alpha p$
\State ${\bf r}_h \gets {\bf r}_h - \alpha q$
\State $z \gets \Call{MultigridCycle}{p_{grid},0,\mathcal{A},{\bf r}_h,\nu_{1},\nu_{2}}$
\State $\beta \gets z^T {\bf r}_h / \delta $
\State $p \gets  {z} + \beta p$ 
\State $\delta \gets {z}^T {\bf r}_h$
\State $i \gets i + 1$
\EndWhile \\
\Return ${\bf u}_h$
\EndFunction
\end{algorithmic}
\end{algorithm}
\FloatBarrier

\subsection{Transfer operators}

The multigrid method relies on the ability to transfer approximate grid solution vectors between nested grids. Transfer operators consists of prolongation operators $\mathcal{P}_c^f:{\Omega}_{c} \rightarrow {\Omega}_f$ that takes a solution vector from the coarse grid to a fine grid, and restriction operators $\mathcal{R}_f^c:{\Omega}_{f} \rightarrow {\Omega}_c$ that takes a solution vector from the fine grid to a coarse grid. The grids are chosen according to $\forall c \in \{1,...,P-1\}, \; \forall f \in \{2,...,P\}, \; c < f$. The natural choice for prolongation and restriction operators in a multigrid method is linear operators.
%
%

In this work we employ the interpolation strategy defined for the first time in \cite{RP87}, where the coarse grid spaces are subspaces of the fine grid space, therefore a solution vector $v_h({\bf x})$ can be expanded in terms of the fine and coarse grid basis functions
\begin{align}
v_h({\bf x})=\sum_{i=1}^{N^{p-1}}v_{i,p-1}\phi_i^{p-1}({\bf x}) = \sum_{j=1}^{N^{p}}v_{j,p}\phi_j^{p}({\bf x}),
\end{align}
and the prolongation as a matrix-vector product interpolating $v$ from $V^{p-1}$ to $V^{p}$
\begin{align}
\mathcal{P}_{p-1}^p v_{p-1} = v_{p}.
\end{align}
Similarly to \cite{RP87}, the restriction operator is a transposition of the prolongation
\begin{equation}
\label{eq:rest}
    {\mathcal{R}}_p^{p-1} = {({\mathcal{P}}^p_{p-1})}^{T}.
\end{equation}





\subsection{Relaxation strategy}

\begin{figure}[H]
\centering \includegraphics[width=0.35\textwidth]{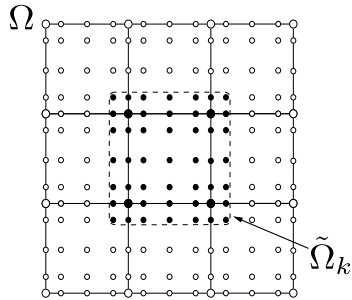}
\caption{An example of one subdomain ${\tilde{\Omega}}_k$ of the Schwarz smoother, with $k$ being the center element in the picture. Subdomain ${\tilde{\Omega}}_k$ is represented by the black nodes, which consist of the element $k$ nodes and the overlapping zone with one overlapping node in each spatial direction. }
\label{fig:schwarzsmoother}
\end{figure}

The second crucial component in a multigrid scheme is the choice of a smoothing operator $\mathcal{S}$, which is responsible for the error correction. Most basic stationary iterative methods, such as Jacobi or Gauss-Seidel schemes prove to work fine for simple test cases on elliptic PDEs (examples in \cite{RP87,FL05}). However, it is known that the performance of the aforementioned methods deteriorates for higher polynomial orders of the finite element discretization \cite{Janssen2011} and for the Galerkin-based frameworks, such as SEM used in this work, efficient smoothers that can exploit element-block structure the aforementioned methods are necessary to ensure rapid decline of the error. One of the most promising strategies for conjugate gradient (CG) discretized schemes is the additive Schwarz method (ASM), which significantly enhances the convergence rate of the multigrid solver in comparison to the standard smoothers \cite{FL05}. The ASM can be presented in the form of a preconditioned defect correction iteration
\begin{align}
    {\bf u}_h = {\bf u}_h - \mathcal{S}^{-1}(\mathcal{A}{\bf u}_h - {\bf f}),
\label{eq:sec_form}
\end{align}
where $\mathcal{S}^{-1}$ is the resulting Schwarz preconditioner. In the context of the multigrid solvers the preconditioner $\mathcal{S}^{-1}$ is a block matrix of the form
\begin{equation}
\label{eq:Mass}
    \mathcal{S}^{-1} = \mathcal{W} \sum_{k=1}^{N^k}{\mathcal{R}^T_{k}(\mathcal{R}_{k}\mathcal{A}\mathcal{R}^T_{k})^{-1}\mathcal{R}_{k}}, 
\end{equation}
where we have introduced a restriction operator $\mathcal{R}_{k}:{\Omega}_{k} \rightarrow {\tilde{\Omega}}_k$ that restricts nodes from the main computational domain ${\Omega}$ to a newly defined subdomain ${\tilde{\Omega}}_k$ (element-block along with the pre-defined overlapping nodes) and $\mathcal{W}$ is diagonal weighting operator
\begin{equation}
    \label{eq:weights}
    \mathcal{W} = \left(\sum_{k=1}^{N^k}{\mathcal{R}_{k}\mathcal{R}^T_{k}}\right)^{-1}.
\end{equation}
In figure \ref{fig:schwarzsmoother} we depict the main Schwarz relaxation strategy used in this work one the exemplary 2D computational domain $\Omega$ consisting of 9 elements. Tests show  that one overlapping node provides a good improvement in convergence rate in comparison to the standard block-Jacobi (ASM with zero overlapping nodes) for a wide range of meshes. In the considered problems in \cref{sec:experiments} we use one overlapping node in each spatial direction for all element-block smoothers as illustrated in figure \ref{fig:schwarzsmoother} unless stated otherwise.

\subsection{Performance of multigrid methods}

To measure the algorithmic efficiency of an iterative multigrid scheme, we can consider the convergence rate. We measure the relative change of the residual of the preconditioned system \eqref{eq:precondsys} on the finest grid between successive iterations
\begin{align}
\label{eq:convrate}
\text{q} = \frac{||{\bf r}_h^{[m]}||_2}{||{\bf r}_h^{[m-1]}||_2}.
\end{align}
The smoothing is done using \eqref{eq:sec_form} and a Work Unit (WU) is measured relative to one sparse matrix vector (SpMV) operation using the system matrix on the fine grid. The multigrid strategy detailed in \cite{ENGSIGKARUP20092100} where semi-coarsening is combined with standard coarsening can be invoked and help improve convergence speed in the context of anisotropic meshes with elements of high aspect ratios.
The assessment of the performance focus on algorithmic efficiency that is assessed in terms of the convergence rate, and the numerical efficiency that is implementation dependent. By selecting the multigrid components it is possible to make trade-offs between algorithmic efficiency and numerical efficiency. 

\subsection{On computational cost of multigrid}

To put a broader perspective on the cost of the solver, we define a relative measurement unit \emph{working unit} (WU). One WU is equivalent to computational cost of one preconditioned defect correction iteration \eqref{eq:sec_form} on the finest grid \cite{Brandt1984}.

We divide the cost of the multigrid method into two parts. First part is the \emph{operator assembly}, it consists of forming a smoothing operator $\mathcal{S}^{-1}$ and coarse grid Laplace operators $\mathcal{A}$ (either by a Galerkin coarse grid operator $\mathcal{RAP}$ or via direct assembly). The second part is the actual solver defined by pseudo-algorithms \ref{alg:mg}, \ref{alg:mgdc} and \ref{alg:mgcg}. The former part is, by a large margin, more expensive. However, in our numerical experiments we pre-compute $\mathcal{S}^{-1}$ and system matrix $\mathcal{A}$ on all grids at $t=0$ s and re-use them throughout the whole time domain, therefore the computational cost of assembling $\mathcal{S}^{-1}$ and $\mathcal{A}$ becomes negligible in compare to the total simulation time.

\section{Numerical benchmarks}
\label{sec:experiments}

We consider two benchmark cases. The first case is a standard test case in 2D for dispersive and nonlinear waves models that has been considered in \cite{EngsigKarup2014} where a high-order finite difference method is used and geometric multigrid method is used in an iterative preconditioned defect correction solver. The second case we consider is a 3D case for one of the blind test experiments corresponding to case 2.3 described in \cite{MaEtAl2015,RansleyEtAl2018}. It represents an advanced case, where the free surface mesh is fully unstructured as is typical for real cases, where a fixed floating production storage and offloading (FPSO) vessel body is represented with curvilinear elements. Remark, different formulations of the Laplace problem is used in each of the benchmarks.

\begin{figure}[H]
\centering
\includegraphics[width=0.9\textwidth]{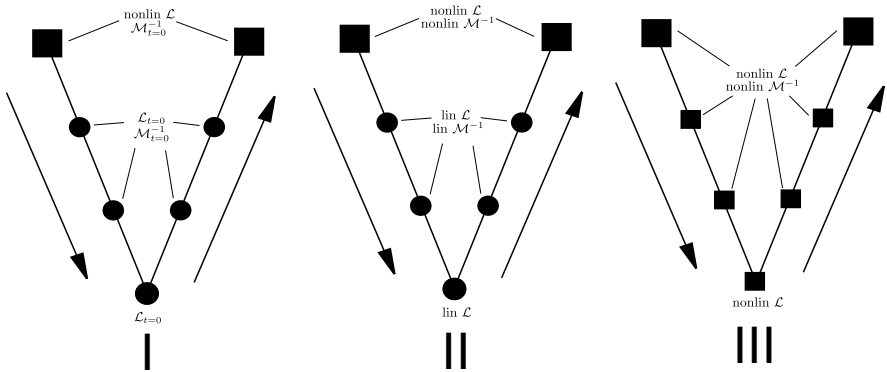}
\caption{Strategies that carry computation throughout the whole time domain. \emph{NONLIN} denotes a full, more computationally expensive, nonlinear Laplace operator \eqref{FSphi} on all grid levels, while the \emph{LIN} is simplified, less expensive, linear version on all levels, but the finest. Subscript '$t=0$' denotes use of an operator from time zero (first step) on all later time steps.} 
\label{fig:simstrats}
\end{figure}

We introduce three different strategies to conduct simulation on the whole time domain as illustrated in figure \ref{fig:simstrats}. To avoid expensive global assembly of the smoothening matrix $\mathcal{S}^{-1}$ and the coarse grid Laplace operators, we store initial $\mathcal{S}^{-1}$ and $\mathcal{A}$ (at $t=0$ s) and apply these in the solver every time step, \emph{I} from figure \ref{fig:simstrats}. We compare the numerical results to the systems \emph{II} and \emph{III}, where we update all the operators on all time steps. Either based on Laplace operator formed from linear formulation of potential flow on coarse grids or based on the formulation \eqref{FSphi} on all grid levels (the latter).  

\subsection{Harmonic generation over a submerged bar}

The test case of harmonic generation over a submerged bar, has become a standard test case for evaluating dispersive and nonlinear wave models. This case was considered in \cite{EEB2016}  where a direct solver for the Laplace problem was used. In this work, we instead consider the geometric $p$-multigrid method and solve the discretized Laplace problem in the $\sigma$-transformed formulation \eqref{eq:laplaceproblemdiscrete}. Qualitatively, we obtain similar results as presented in the previous work, hence, we report the results obtained from using the iterative multigrid solver strategies.

We consider a numerical setup where the time domain is $T:t=(0,82.8)$ s resolved with a fixed time step $\Delta t=0.0736$ s and the spatial fluid domain with high aspect ratio elements  $\Omega = (0,29) \; \text{x} \; (0,1)\; $ m$^2$. The mesh consists of quadrilateral elements with number of elements in each direction such that  $N^k_x=103$, $N^k_y=1$, $P_x = P_y = 6$. The discretization of the linear system results in discrete anisotropy of $\frac{dz}{dx} \sim 8$ where $dx$ and $dz$ refers to the element size in each of the Cartesian directions. Natural restriction strategy would be to semi-coarse in $x$ direction first to perform more smoothing operations on better conditioned system. In this case however, the assembled block-smoother has proven to be efficient enough. The coarsening strategy is given in table \ref{table:BTcoarse}. 
Results were computed using a fourth-order explicit Runge-Kutta (ERK4) method for the temporal integration. In total the simulation run for 1125 time steps and the results are averaged on a total of Laplace solves of 5625. Two pre- and post-smoothening are used for the condensed system, and one pre- and post smoothening for the full system of equations.

\begin{figure}[H]
\subfloat[ \label{subfig:t_scaling}]{%
  \includegraphics[width=0.45\textwidth]{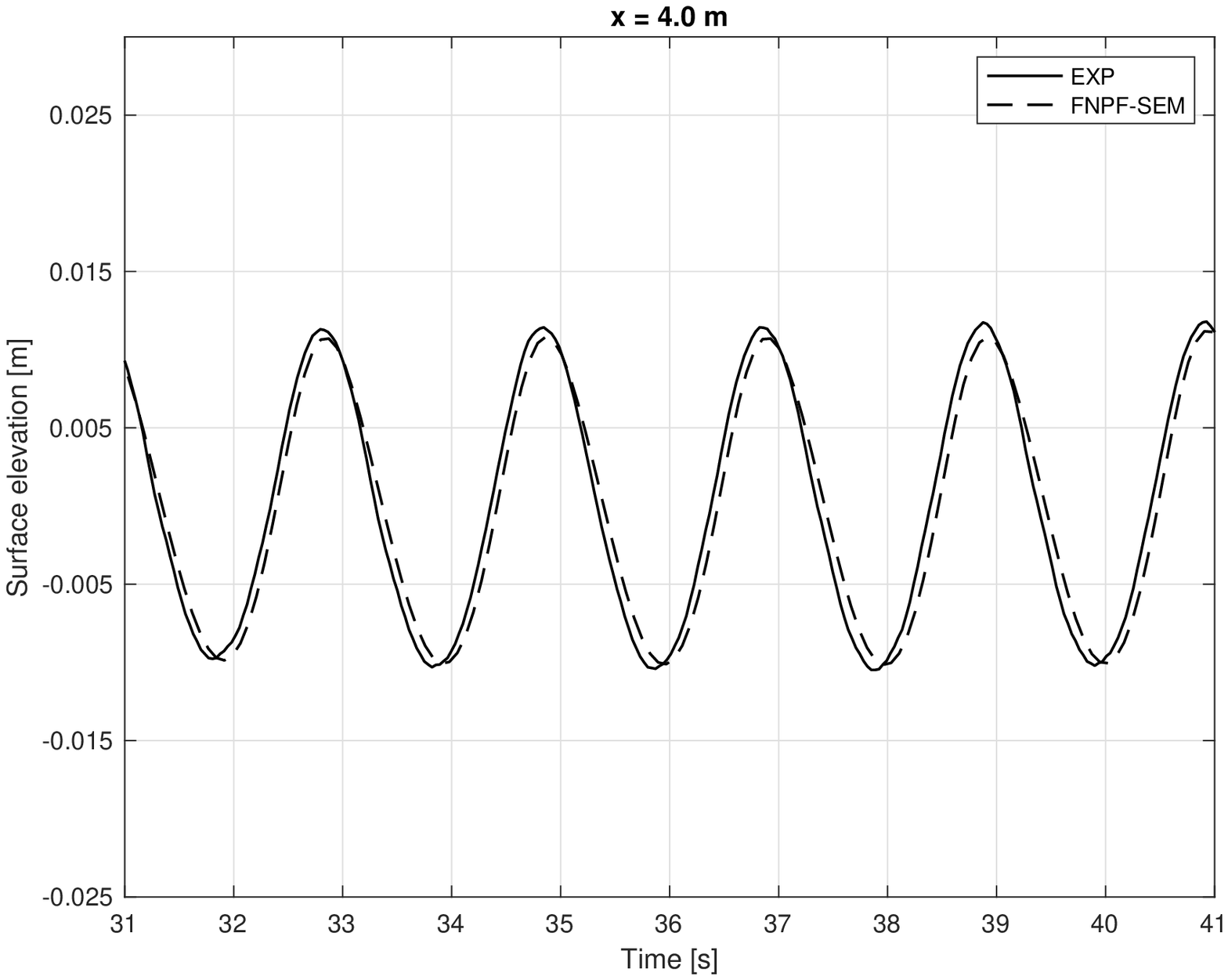}
}
\hfill
\subfloat[\label{subfig:q_scaling}]{%
  \includegraphics[width=0.45\textwidth]{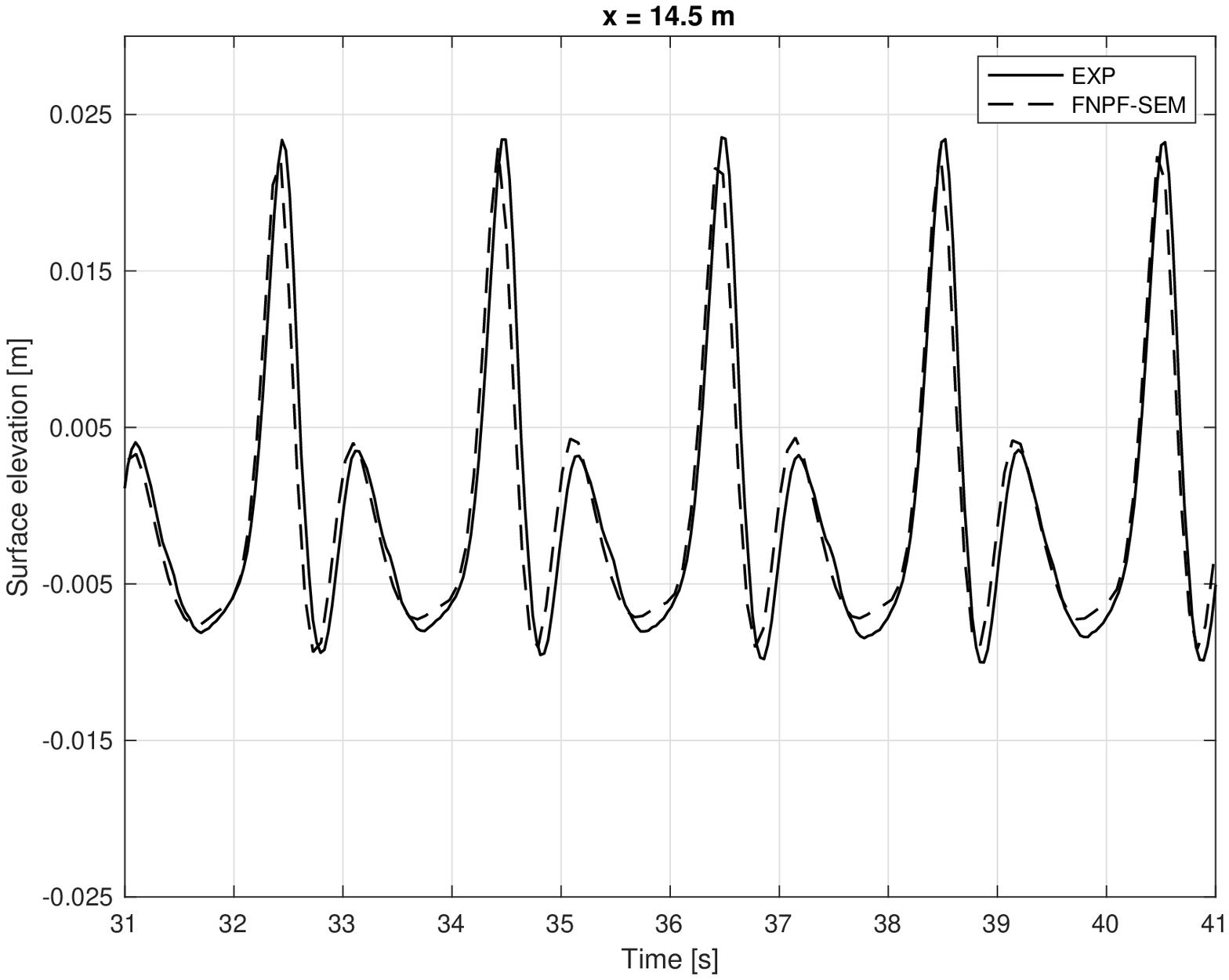}
}
\vfill
\subfloat[ \label{subfig:t_scaling2}]{%
  \includegraphics[width=0.45\textwidth]{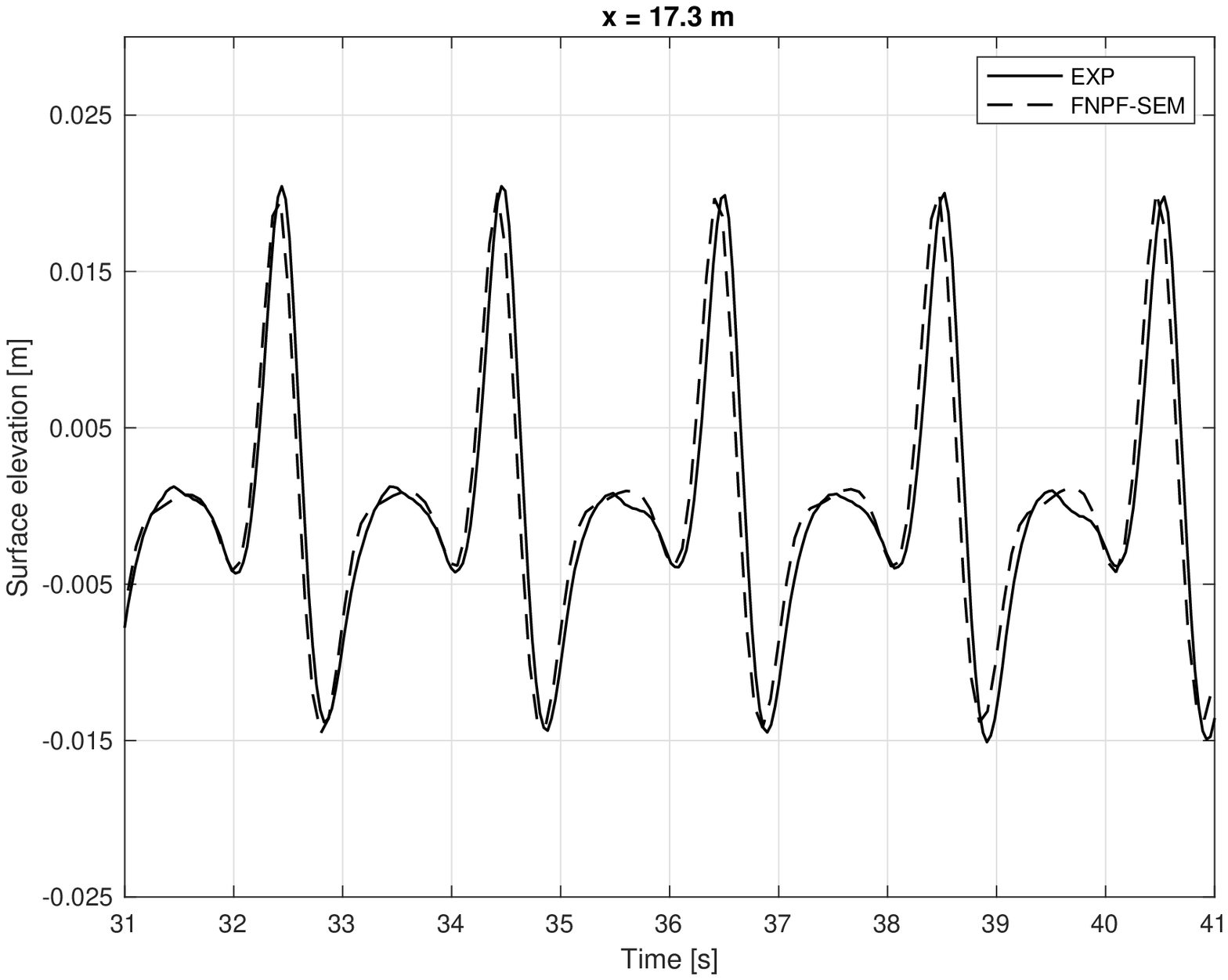}
}
\hfill
\subfloat[\label{subfig:q_scaling2}]{%
  \includegraphics[width=0.45\textwidth]{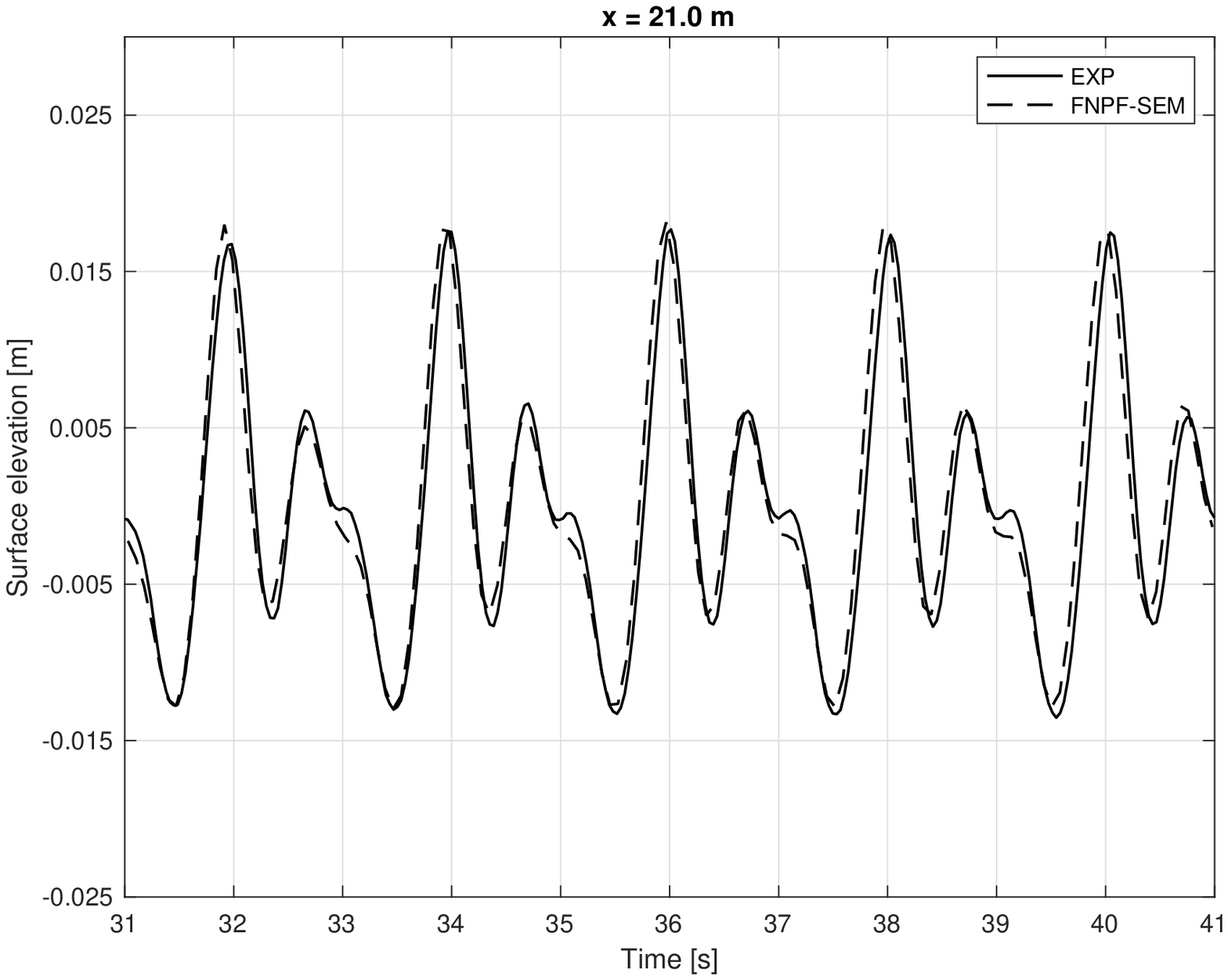}
}
\caption{Time series histories of free surface elevations  in comparison with experimental measurements data at wave gauges positioned at $x$=4.0m, 14.5m, 17.3m and 21.0m.}
\label{fig:smoothing}
\end{figure}

\FloatBarrier
\begin{table}[!htbp]
\centering
\caption{Coarsening strategy.}
\label{table:BTcoarse}
\begin{tabular}{@{}lccc@{}} \toprule
Grid level & 3 & 2 & 1 \\ \midrule
 $P_x$ & 6 & 3 & 1 \\
 $P_y$ & 6 & 3 & 1 \\
 \bottomrule
\end{tabular}
\end{table}
\FloatBarrier

Tables \ref{table:BarTest_iterf} and \ref{table:BarTest_WUf} depict results for several different schemes, that can be split with in terms of the applied iterative method: multigrid preconditioned conjugate gradient (PCG-$p$-GMG)
, multigrid preconditioned defect correction (PDC-$p$-GMG)
, stand-alone multigrid ($MG$) and also to the unique strategies presented in figure \ref{fig:simstrats}. Table \ref{table:BarTest_iterf} depicts average iteration number and table \ref{table:BarTest_WUf} the measured time (in WUs) to achieve a given tolerance. The results are compared to a direct solver based on Cholesky factorization of the system matrix with permutation of both rows and columns using a reverse Cuthill McKee (RCM) ordering.

\FloatBarrier
\begin{table}[!htbp]
\centering
\caption{Average number of iterations to achieve given tolerance. Full system.}
\label{table:BarTest_iterf}
\begin{tabular}{@{}lcccc@{}} \toprule
Method & $1e-4$ & $1e-5$ & $1e-6$ & $1e-7$ \\ \midrule
PCG-$p$-GMG(3,3) I & 1.00 & 1.99 & 1.99 & 2.00 \\
PCG-$p$-GMG(3,3) II & 1.00 & 1.99 & 1.99 & 2.00 \\
PCG-$p$-GMG(3,3) III & 1.00 & 1.99 & 1.99 & 2.00 \\
MG(3,3) I & 1.00 & 1.99 & 1.99 & 2.01 \\
PDC-$p$-GMG(3,3) I & 1.00 & 1.99 & 1.99 & 2.01 \\
\bottomrule
\end{tabular}
\end{table}
\FloatBarrier

\FloatBarrier
\begin{table}[!htbp]
\centering
\caption{Average cost to achieve given tolerance in absolute WUs $[ms]$.}
\label{table:BarTest_WUf}
\begin{tabular}{@{}lcccc@{}} \toprule
Method & $1e-4$ & $1e-5$ & $1e-6$ & $1e-7$ \\ \midrule
PCG-$p$-GMG(3,3) I & 8.24 & 16.52 & 15.83 & 17.04 \\
PCG-$p$-GMG(3,3) II & 8.44 & 16.71 & 16.21 & 16.62 \\
PCG-$p$-GMG(3,3) III & 8.66 & 16.68  & 16.30 & 16.36 \\
MG(3,3) I & 8.67 & 16.93 & 16.84 & 16.86 \\
PDC-$p$-GMG(3,3) I & 8.25 & 15.71 & 16.97 & 17.38 \\
Direct, LU & \multicolumn{4}{c}{73.75}  \\
\bottomrule
\end{tabular}
\end{table}
\FloatBarrier

\begin{figure}[!htbp]
\centering
\subfloat[\label{subfig:WG7}]{%
\includegraphics[width=0.9\textwidth]{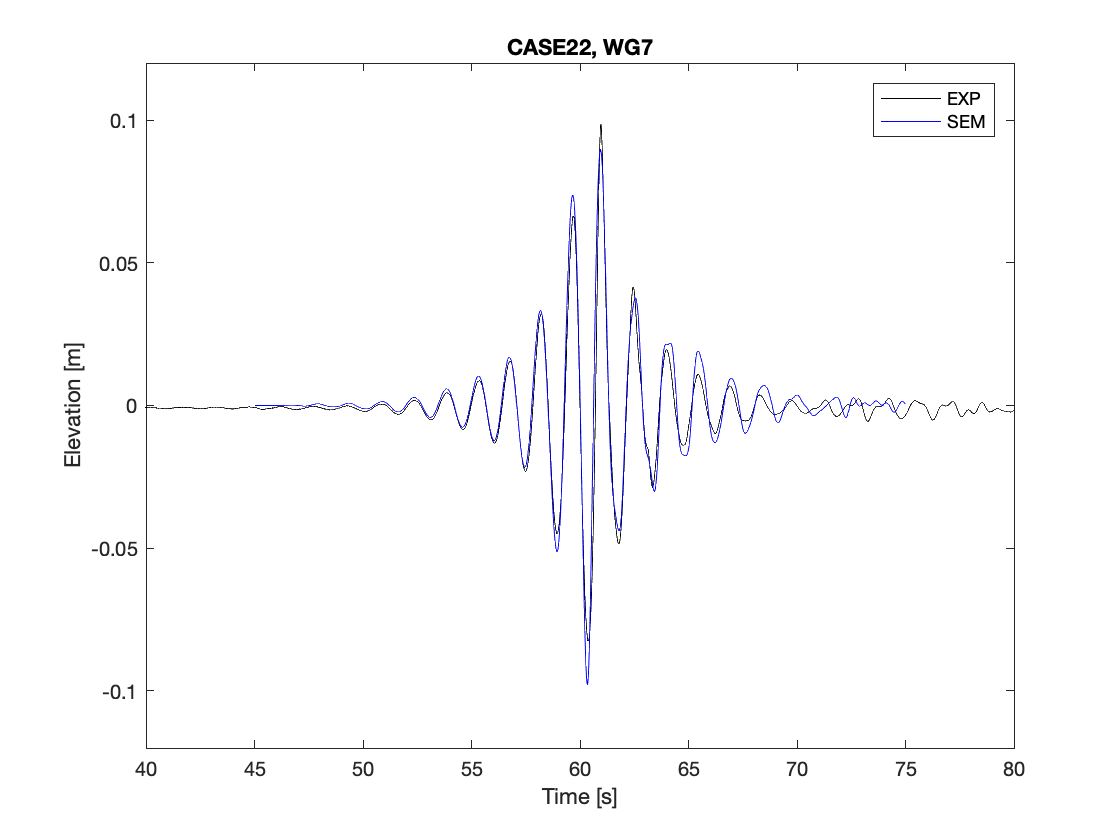}
}
\vfill
\subfloat[ \label{subfig:WG24}]{%
\includegraphics[width=0.9\textwidth]{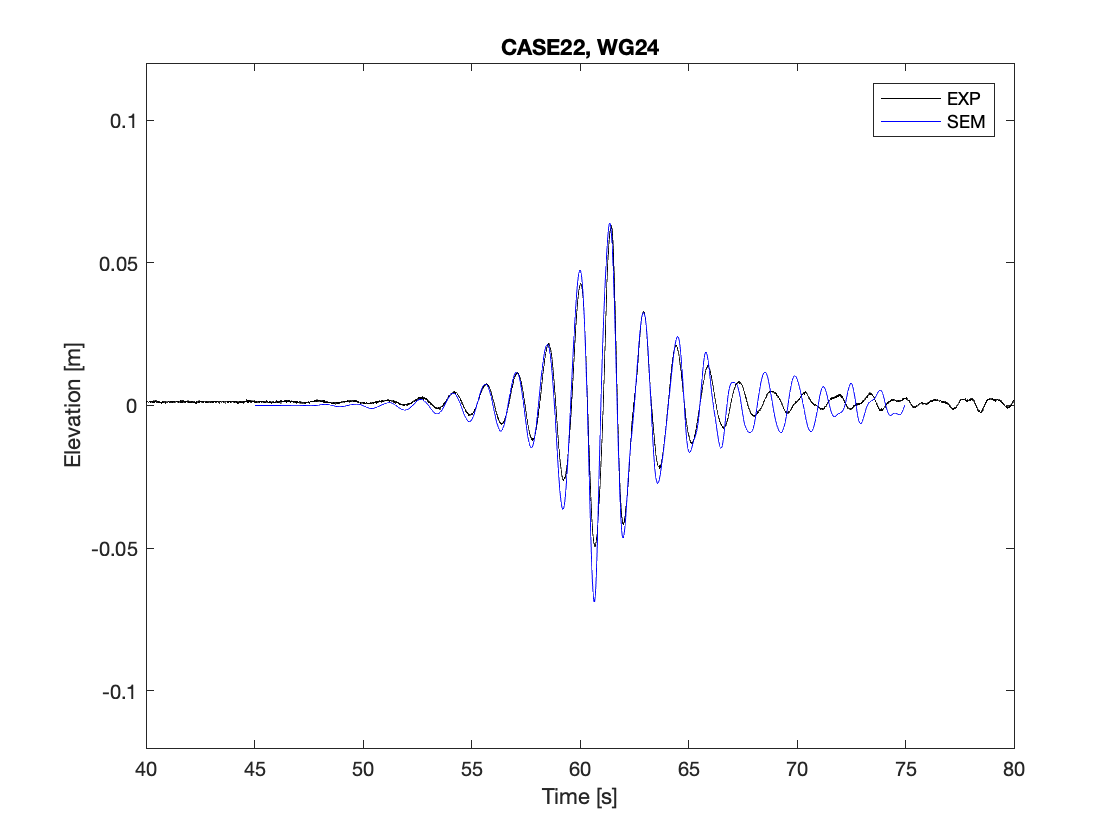}
}
\caption{Wave gauges measurements for a focusing wave interacting with a fixed FPSO. A comparison is made with experimental results from the  CCP-WSI case 2.2 for wave gauges a) WG7 and b) WG24.}
\label{fig:WPmeasurements}
\end{figure}

The $p$-GMG strategies gives similar results in solver performances, and is demonstrated in 2D to give speedups ranging from 4.5 to 8.9 times depending on choice of tolerance levels. Although for NONLIN cases it does not differ much from other  strategies, when a Laplace operator based on assuming $\eta=0$ (small-amplitude wave assumption) is deployed (worse but cheaper approximation), it is the only strategy among LIN that maintains TME. More universal way of representing the results in table \ref{table:BarTest_iterf} is in Working Units (WU), cf. table \ref{table:BarTest_WUf}. One WU is set to be a cost of a sparse matrix vector (SpMV) operation. 
Figures \ref{fig:WPmeasurements} present exemplary results of the computed free surface time series compared with experimental values released in connection with the blind test experiment described in \cite{RansleyEtAl2018}. 

Remark, the number of outer CG iterations can be compared to results given in \cite{CaiEtAl1998} and shows to be more efficient than the reported FEM results, with close to an order of magnitude improvement in iteration counts.

%


\subsection{Focusing wave impacting a FPSO structure in 3D}
\label{sec:fpso3d}

We consider numerical experiment where we compare to the case 2.3 of the experimental measurement campaign from the FROTH (EPSRC) project \cite{MAI2016115} carried out at Plymouth University in the COAST Laboratory facility in UK. In case 2.2 the nonlinear focusing wave group is determined using New Wave Theory \cite{TromansEtAl1991} and is based on a JONSWAP spectrum with parameters $A=0.08930$ m, $T_p=1.456$ s, $h=2.93$ m, $Hs=0.103$ m, $kA=0.17$ and  $\alpha=0.174533$ rad (angle of propagation relative to the centre-line of the basin). The measurement data CCP-WSI ID 22BT1 is provided via the CCP-WSI project \cite{CCPWSI1}. The wave signal is generated by reverse engineering using harmonic analysis, and we use the procedure based on FFT for reproducing the time series of WG1 as described \cite{TromansEtAl1991}.

\begin{figure}[!htbp]
\centering
\subfloat[\label{subfig:FPSO1}]{%
  \includegraphics[width=0.95\textwidth]{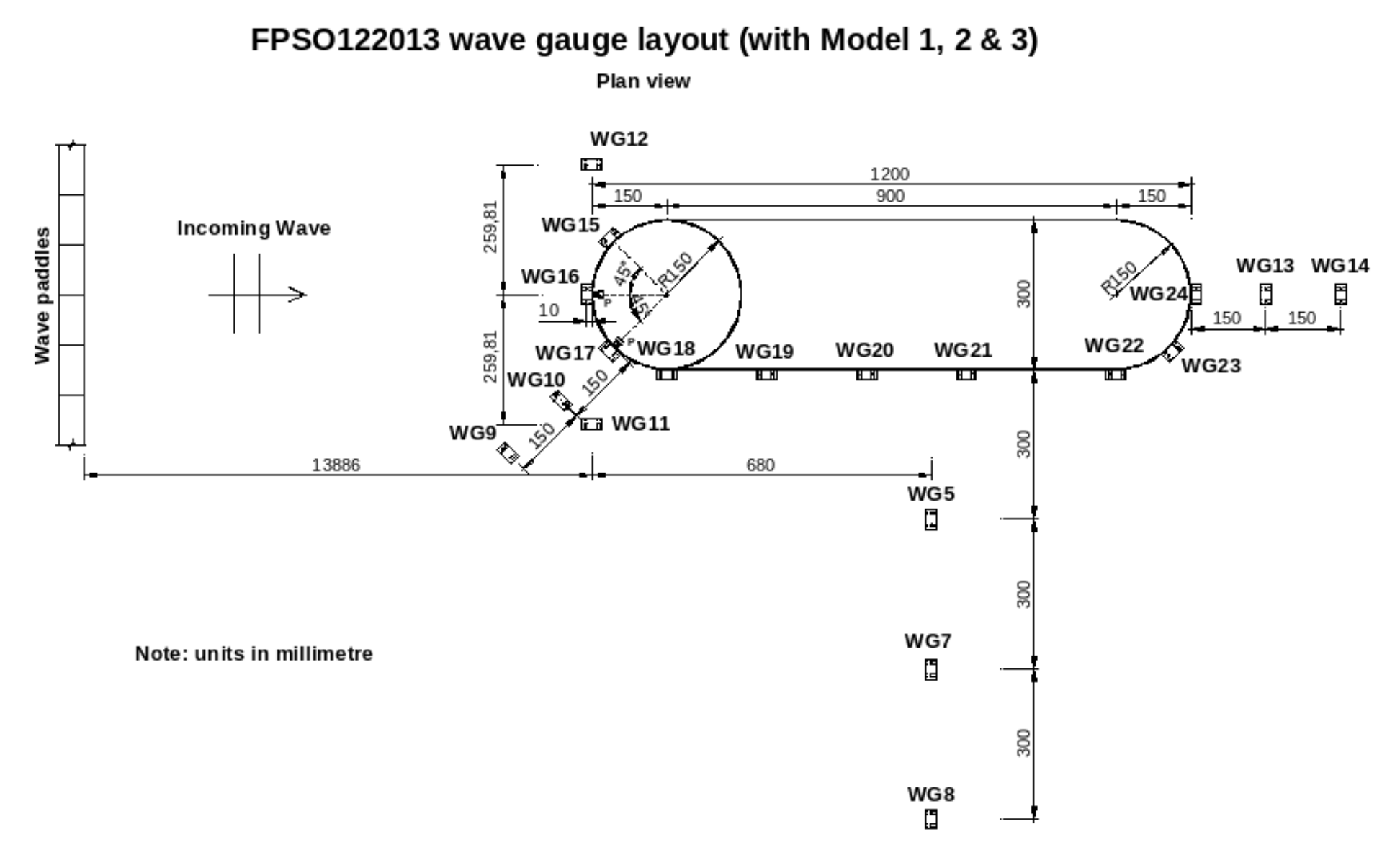}
}
\vfill
\subfloat[ \label{subfig:FPSO2}]{%
\includegraphics[width=0.95\textwidth]{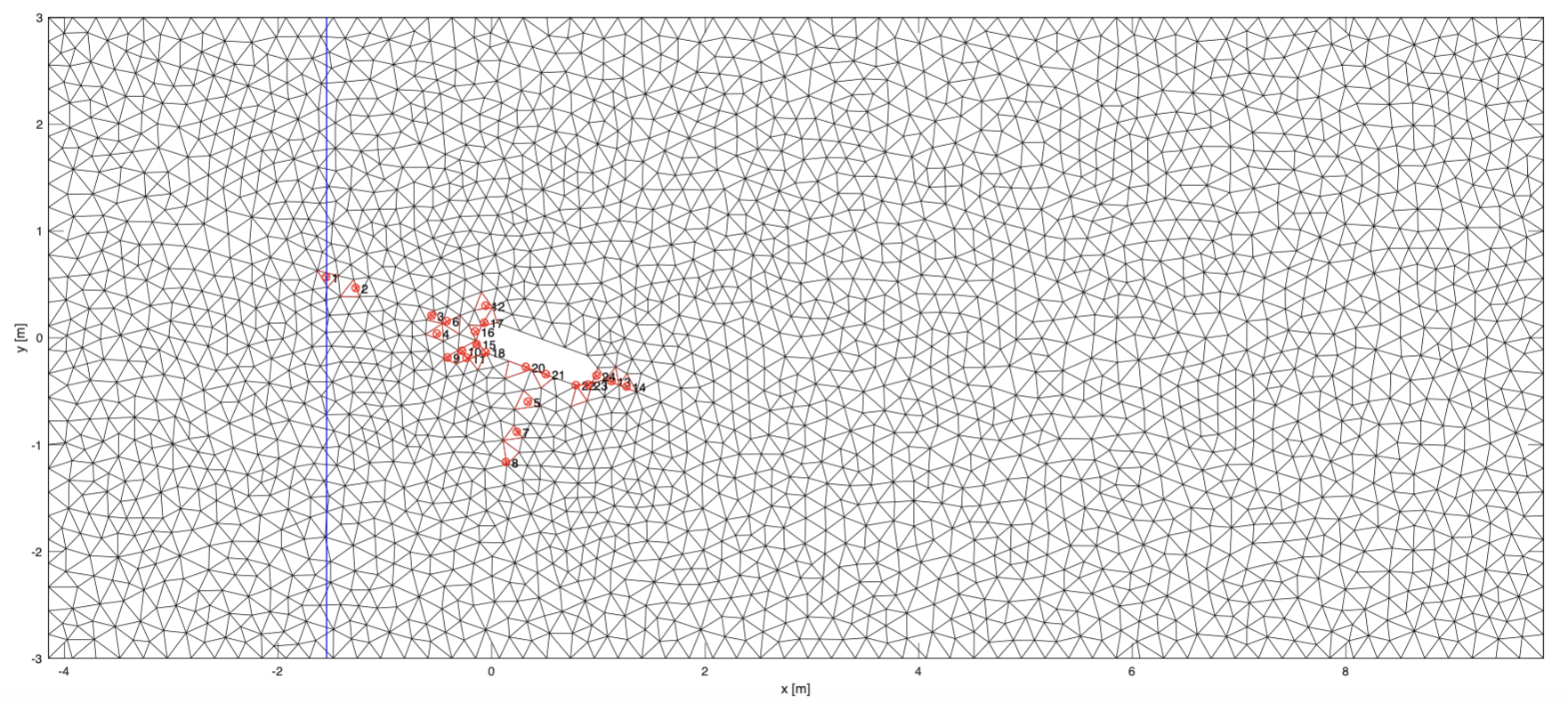}
}
\caption{The problem setups illustrated for the nonlinear focusing wave group interacting with a fixed FPSO at an angle. a) Experimental setup with wave generation and wave gauges. b)  Unstructured mesh setup for the FPSO at an incident angle of 20 degrees with signal measurement wave gauges indicated. The (blue) interface is positioned at at wave gauge 1 (WG1) and is highlighted between the relaxation zone (left of interface) and the computational domain (right of interface).}
\label{fig:setup}
\end{figure}

We employ a FNPF-SEM model that is based on a non-$\sigma$-transformed formulation that results from the SEM Galerkin discretization of \eqref{eq:laplaceeq}. The time series signal of the nonlinear focusing wave group is compared to the experimental data from different wave gauges in the time domain $T:t\in[40,80]$ s. In the simulation, a fixed time step $\Delta t=0.025$ s is used and the spatial domain is defined with a high aspect ratio of elements for a numerical wave tank of size $\Omega = [-14.036,18.364]\times[-0.125,0.125]\times[-2.93,0]$ m$^3$ with a fixed FPSO embedded in the fluid domain. The mesh setup consist $N^k = 864$ 3D prism elements (see figure \ref{fig:setup}) and polynomial orders $P_{xy} = 5$ (triangular base) and $P_z = 3$. Semi-coarsening in the multigrid solver is the horizontal directions is employed until polynomial orders in all directions are equal in the expansion basis. Then standard coarsening is employed until the coarsest level is reached. This is illustrated in table \ref{table:BTcoarse3D}. On all $p$-GMG levels, but the coarsest grid, three pre- and post-smoothening operations are used. Results were computed using the ERK4 for the temporal integration. In total the simulation run 1\;200 time steps. 

\FloatBarrier
\begin{table}[!htbp]
\centering
\caption{Coarsening strategy.}
\label{table:BTcoarse3D}
\begin{tabular}{@{}lccc@{}} \toprule
Grid level & 3 & 2 & 1 \\ \midrule
 $P_{xy}$ & 5 & 3 & 1  \\
 $P_z$ & 3 & 3 & 1 \\
 \bottomrule
\end{tabular}
\end{table}
\FloatBarrier

Similarly to the previous submerged bar test case, tables \ref{table:FPSO3D_iterf} and \ref{table:FPSO3D_WUf} present the numerical results of the PCG-$p$-GMG, PDC-$p$-GMG and $MG$ solvers. Unlike the 2D problem, it can be seen in figure \ref{table:FPSO3D_iterf} that for the 3D case number of iterations deteriorates with lower tolerance levels, however, all the solvers still maintain approximately 1 order of magnitude per iteration convergence rate. The results in table \ref{table:FPSO3D_WUf} demonstrates speedup in comparison to a sparse direct solver with reordering ranging from 0.9 (slowdown) to 7.3 times depending on the choice of tolerance level. The slowdown at strict tolerance is expected as the benefits of using multigrid reduces with increasing number of iterations. However, we note that in contrast to a sparse direct solver that for highly unstructured meshes sees an impact on efficiency, the iterative $p$-multigrid solver strategy is designed to achieve $\mathcal{O}(n)$ complexity. Therefore, the results of this test case confirms the feasibility and potential of using an iterative $p$-multigrid spectrel element model on an advanced application to improve efficiency.   

\FloatBarrier
\begin{table}[!htbp]
\centering
\caption{Average number of iterations to achieve given tolerance.}
\label{table:FPSO3D_iterf}
\begin{tabular}{@{}lcccc@{}} \toprule
{\bf Method} & $1e-4$ & $1e-5$ & $1e-6$ & $1e-7$ \\ \midrule
PCG-$p$-GMG(3,3) I   & 1.11 & 2.33 & 5.06 & 8.47\\
PCG-$p$-GMG(3,3) II  & 1.11 & 2.33 & 5.08 & 8.51\\
PCG-$p$-GMG(3,3) III & 1.11 & 2.33 & 5.09 & 8.52\\
$p$-GMG   (3,3) I   & 1.11 & 2.33 & 5.08 & 8.69\\
PDC-$p$-GMG(3,3) I & 1.11 & 2.33 & 5.09 & 8.69\\
\bottomrule
\end{tabular}
\end{table}
\FloatBarrier

\FloatBarrier
\begin{table}[!htbp]
\centering
\caption{Average cost to achieve given tolerance in relative WUs $[ms]$.}
\label{table:FPSO3D_WUf}
\begin{tabular}{@{}lcccc@{}} \toprule
{\bf Method} & $1e-4$ & $1e-5$ & $1e-6$ & $1e-7$ \\ \midrule
PCG-$p$-GMG(3,3) I  	& 8.09 & 17.03 & 40.62 & 64.84\\
PCG-$p$-GMG(3,3) II 	& 8.35 & 17.24 & 41.25 & 65.11\\
PCG-$p$-GMG(3,3) III	& 8.73 & 18.01 & 42.35 & 67.20\\
$p$-GMG   (3,3) I    & 8.16 & 16.89 & 41.49 & 65.00\\
PDC-$p$-GMG(3,3) I    & 8.01 & 16.59 & 41.35 & 66.09\\
Direct, LU              & \multicolumn{4}{c}{59.11}  \\
\bottomrule
\end{tabular}
\end{table}
\FloatBarrier

\subsection{Scaling of the geometric p-multigrid}
\label{sec:scaling}

The scalability properties of the solver are examined in terms of increasing number of elements $N^k$ and polynomial order $P$. We consider a deep water nonlinear standing wave problem using the approximate 2D solution due to Glozmann \cite{AG96} that is used in a 3D tank and only serve to run this benchmark using a nonlinear wave problem in a finite domain. The wave height is $H=0.089$ m, wave length is $L=2$ m and the dispersion parameter is $kh=2\pi$, where we use a flat bottom with still water depth $h=2$ m on the spatial domain $\Omega = (0,8) \; \text{x} \; (-4,4) \; \text{x} \; (-1,0)\; $ m$^3$ with no body. 

To investigate scalable properties of the solver with mutable number of elements $N^k$ we employ a series of five 3D prismatic $h$-meshes with the following number of elements on each $N^k = 40,252,1134,2712,5088$ with approximately constant aspect ratio $\frac{d_{xy}}{d_z} \approx 3$. Elements in each mesh are further discretized by the same pair of polynomial orders $P_{xy}=5$ and $P_z=3$. The geometric $p$-multigrid follows the  same multigrid strategy as described in Section \ref{sec:fpso3d}. To test the scalability property, we consider a set of polynomial order pairs ($P_{xy},P_z)=(2,2),(4,2),(5,3),(6,4),(9,7)$ that discretize the problem on the finest grid. We employ the same coarsening strategy for all cases. First, the problem is semi-coarsened in the vertical direction to reach the polynomial order $P_{xy} = P_z = P$, and thereafter standard coarsening is used according to $P_{coarse} = {ceiling}((P+1)/2)$ down to $P=1$.

\begin{figure}[!htbp]
\centering
\includegraphics[width=0.45\textwidth]{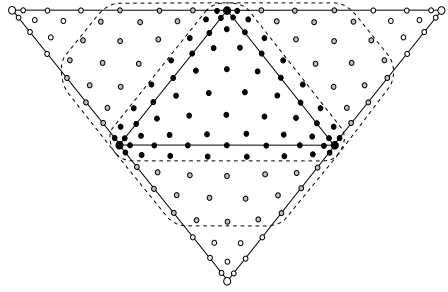}
\caption{Exemplary overlapping regions on Fekete points for polynomial orders of the expansion on elements $P_x=P_y=8$. Filled black circles represent overlapping region of one node in all spatial directions (constant for all polynomial orders in figure \ref{fig:smoothing}. Black with addition of filled grey circles represent \textit{refined Schwarz} smoother with  overlapping nodes in all spatial directions.}
\label{fig:feketeschwarz}
\end{figure}

\begin{figure}[!htpb]
\subfloat[ \label{subfig:t_scaling3}]{%
  \includegraphics[width=0.48\textwidth]{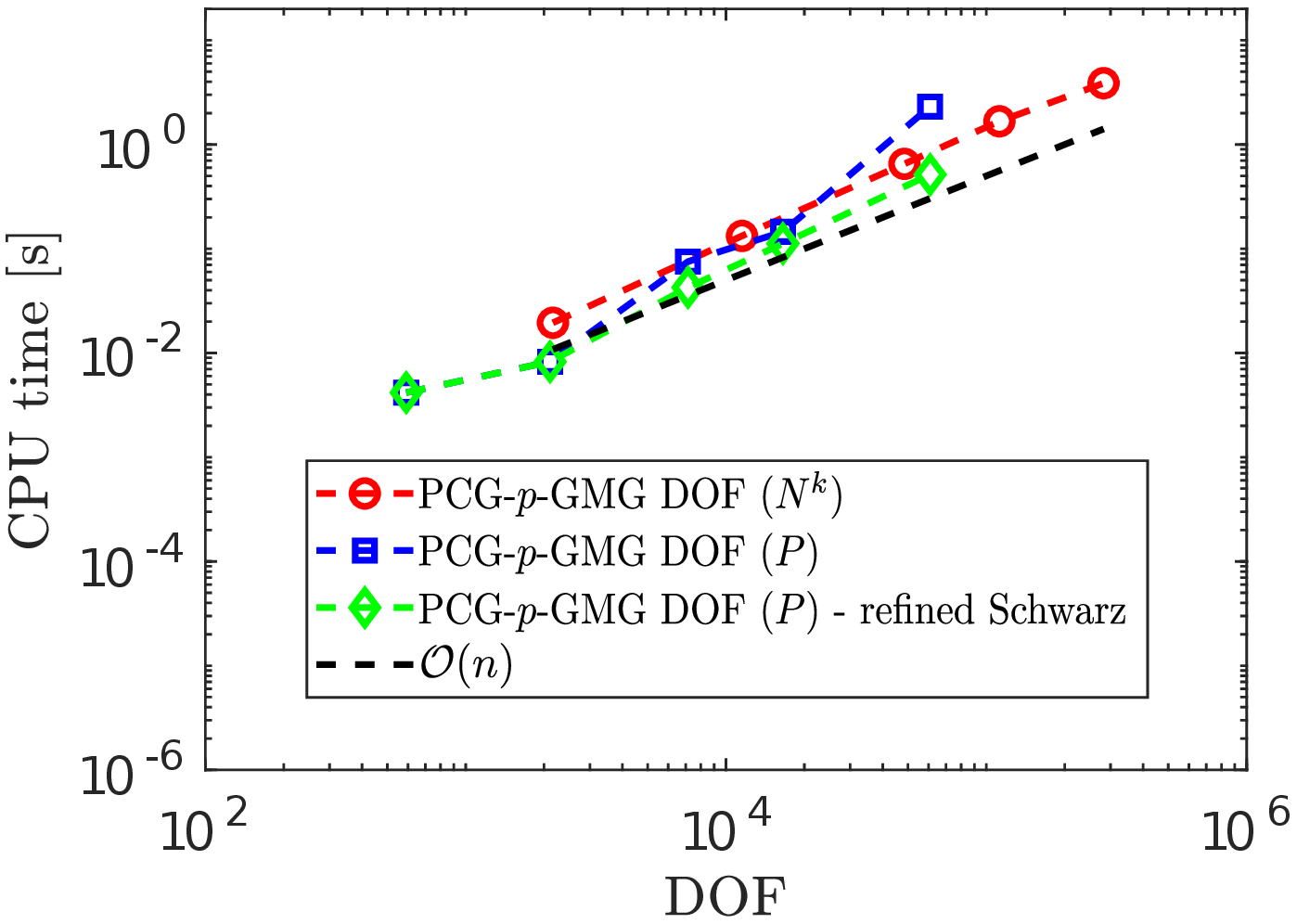}
  }
\hfill
\subfloat[\label{subfig:q_scaling3}]{%
  \includegraphics[width=0.48\textwidth]{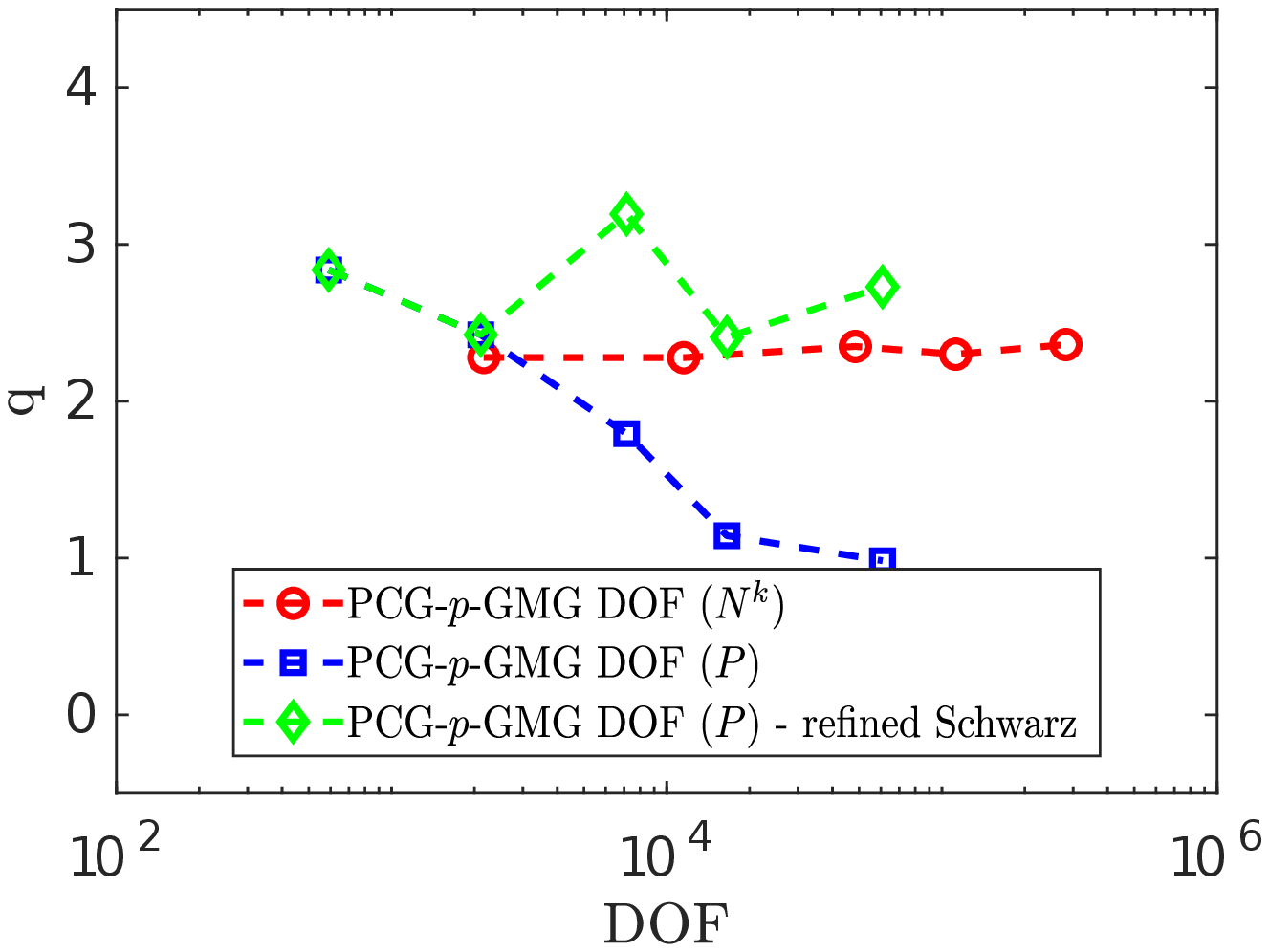}
}
\caption{Scaling of computational time and convergence rate of the geometric $p$-multigrid solver with respect to degrees of freedom (DOF) of the linear system.}
\label{fig:scaling}
\end{figure}

Figure \ref{fig:scaling} reports computational times and convergence rates for all the considered problems with various $N^k$ and $P$. Convergence rate $\text{q}$ is defined by equation \eqref{eq:convrate}. To test the scalability property, we consider two relaxation strategies: first with constant number of overlapping nodes 1 for all cases and all coarse grids and second one (named \textit{refined Schwarz} in figure \ref{subfig:q_scaling}) where the number of overlapping nodes is bounded by ${ceiling}((P+1)/2)$ with $P$ being the polynomial order on the respective grid. Figure \ref{fig:feketeschwarz} presents both strategies on Fekete nodes considered in the 3D prismatic mesh. 

As seen in figure \ref{fig:scaling}, the geometric $p$-multigrid is shown to have excellent scalability properties with the number of elements $N^k$, where the convergence rate $\text{q}$ stays constant on all considered meshes, which results in perfectly linear computational cost increase, cf. figure \ref{subfig:t_scaling}. A good convergence rate is achieved on various cases for polynomial orders $P$ with the refined Schwarz method. The convergence measure $\text{q}$ given in \eqref{eq:convrate}, however, shows a some variation in the various setups subject to varying number of elements $N
^k$ in contrast to those with fixed polynomial order $P$. This comes with increased cost of the preconditioner matrix-vector operations. Despite this, the solver achieves  $\mathcal{O}(n)$ (with $n$ referring to the degrees of freedom) for the scaling of the computational cost. The solver with non-refined Schwarz smoother has a visible convergence rate loss, which results in increasing computational cost.  


\FloatBarrier
\FloatBarrier


\section{Conclusions}
\label{sec:conclusions}

We propose and demonstrate efficient and scalable iterative solution of the Laplace equation that arises in FNPF models. We propose to use a $p$-type geometric multigrid ($p$-GMG) iterative solver as the preferred way to solve the FNPF model equations, since it provides a natural and efficient way of refining the grid resolution and maintaining as much information as possible of the geometric features of the structural body representation across all nested grids in the multigrid solver.

We have tested several geometric $p$-multigrid accelerated iterative strategies for efficient solution of the Laplace problem discretized using a high-order spectral element method. The main objective was to deliver textbook multigrid efficiency (TME) for spatial discretization using the spectral element method on general unstructured meshes to be able to handle wave-body problems. In the current implementations, we have shown that it is possible to achieve TME at solver tolerances around $10^{-4}$ and less than two times TME for tolerances down to $10^{-7}$ that would be relevant to reduce the algebraic errors throughout a simulation. In numerical benchmarks of practical relevance, we have demonstrated for a FNPF-SEM solver that it is possible to achieve engineering accuracy (within tolerances of \textrm{atol} = $10^{-5}$ and \textrm{rtol} = $10^{-4}$). This leads to a main practical result of using $p$-GMG solvers, namely, that it is possible to achieve a low average iteration count of at most 2 iterations for the stage solver in explicit Runge-Kutta methods for all tolerances larger than $10^{-7}$. We tested both FNPF modelling using a traditional $\sigma$-transformed formulation suitable for pure wave propagation as well as non-$\sigma$-formulation suitable for wave-body applications. 

Based on the presented numerical experiments, the speedup in the case of a nonlinear focusing wave group presented up to 7.3 times speedup using a spectral $p$-multigrid solver over a sparse direct solver (cf. section \ref{sec:fpso3d}). The resulting FNPF-SEM simulator is as a main result efficient and prepared for advanced simulations such as a wave-body problem with an FPSO subject to wave-induced loads. The proposed iterative solver strategy demonstrated for a PCG-$p$-GMG iterative solver linearly scalable work effort (cf. section \ref{sec:scaling}), which is the key requisite for enabling large-scale practical calculations. In ongoing work, the iterative $p$-multigrid method is implemented in a parallel computing framework via domain decomposition designed for massive scalability to prepare it for execution on massively parallel computing systems.





%

\section*{Acknowledgments}

The research was hosted by Department of Applied Mathematics and Computer Science (DTU Compute) and access to high-performance computing resources granted via the DTU Computing Center. 

\bibliographystyle{siamplain}
\bibliography{references}
\end{document}